\documentclass[11pt]{article}
\usepackage[utf8]{inputenc}
\usepackage{amsmath,amssymb,amsthm,amsfonts}
\usepackage{enumitem}
\usepackage{geometry}
\usepackage{authblk}
\usepackage[hidelinks]{hyperref}
\setlist[enumerate]{font=\upshape}
\numberwithin{equation}{section}

\theoremstyle{plain}
\newtheorem{theorem}{Theorem}[section]
\newtheorem{lemma}[theorem]{Lemma}
\newtheorem{proposition}[theorem]{Proposition}
\newtheorem{corollary}[theorem]{Corollary}
\theoremstyle{definition}
\newtheorem{remark}[theorem]{Remark}

\newcommand{\Q}{\mathbb Q}
\newcommand{\Fq}{\mathbb F}
\newcommand{\alg}{\mathrm{alg}}
\DeclareMathOperator{\Aut}{Aut}
\DeclareMathOperator{\Gal}{Gal}
\DeclareMathOperator{\Br}{Br}
\DeclareMathOperator{\inv}{inv}
\DeclareMathOperator{\ind}{ind}
\DeclareMathOperator{\tr}{tr}
\DeclareMathOperator{\Trd}{Trd}
\DeclareMathOperator{\Nrd}{Nrd}
\DeclareMathOperator{\GL}{GL}
\DeclareMathOperator{\SL}{SL}
\DeclareMathOperator{\SK}{SK}
\DeclareMathOperator{\End}{\mathsf{End}}

\title{An example of a non-commutative division ring with abelian maximal subgroups}
\author[1]{Huynh Viet Khanh}
\affil[1]{Department of Mathematics and Informatics, Ho Chi Minh City University of Education, Ho Chi Minh City, Vietnam\\\texttt{khanhhv@hcmue.edu.vn}}
\date{}
\hypersetup{
  pdftitle={An example of a non-commutative division ring with abelian maximal subgroups},
  pdfauthor={Huynh Viet Khanh}
}

\begin{document}
\maketitle

\begin{abstract}
We construct a cyclic division algebra $D$ of degree $6$ with center $F$ algebraic over $\Q$. The multiplicative group $K^\times$ of a non-Galois maximal subfield $K$ is maximal in $D^\times$, disproving Conjecture~C of Akbari, Ebrahimian, Momenaee Kermani, and Salehi Golsefidy. The field norm $N_{K/F}$ maps $K^\times$ onto $F^\times$, and its kernel is maximal in the reduced norm-one group $D^{(1)}$. The quotient $D^{(1)}/Z(D^{(1)})$ is nonabelian simple, with the image of $\ker N_{K/F}$ as an abelian maximal subgroup. The proof uses matrix factorizations over local fields, results of Breuillard and Glazer on fibers of convolution powers, and Moret--Bailly's theorem on points over extensions with prescribed splitting.
\end{abstract}

{\small
\begin{quote}

\end{quote}
}

\section{Introduction}

Akbari, Mahdavi-Hezavehi, and Mahmudi \cite[Conjecture~2, p.~431]{AMM1999} conjectured that a division ring $D$ is commutative if its multiplicative group $D^\times$ has a nilpotent maximal subgroup. This conjecture and its abelian case are stated as Conjectures~B and~C, respectively, in Akbari et al.\ \cite[p.~202]{AEMSG}. Thus, Conjecture~C asserts that if $D^\times$ has an abelian maximal subgroup, then $D$ is commutative. In this note, we construct a cyclic division algebra of degree~$6$ whose multiplicative group has abelian maximal subgroups. These subgroups are nilpotent, so the example disproves both Conjectures~B and~C.

In order to state the main result, we recall some notation. A division ring is said to be \textit{centrally finite} if it is finite dimensional over its center, as in Lam \cite[\S14, Def.~(14.1), p.~217]{Lam2001}. Let $D$ be a centrally finite division ring with center $F$, and let $K$ be a maximal subfield of $D$. We write $N_{K/F}\colon K^\times\to F^\times$ and $\Nrd_D\colon D^\times\to F^\times$ for the field norm and the reduced norm, respectively, and we let
\[
D^{(1)}=\ker\Nrd_D=\{d\in D^\times\mid\Nrd_D(d)=1\}.
\]
Since $\Nrd_D|_{K^\times}=N_{K/F}$, we have
\begin{equation}\label{eq:kerN}
\ker N_{K/F}=\{k\in K^\times\mid N_{K/F}(k)=1\}=K^\times\cap D^{(1)}.
\end{equation}
Thus, $\ker N_{K/F}\trianglelefteq K^\times\leq D^\times$ and $\ker N_{K/F}\leq D^{(1)}\trianglelefteq D^\times$. We use the notation $\Nrd_A$ and $A^{(1)}$ in the same sense for any central simple algebra~$A$.

\begin{samepage}
\begin{theorem}\label{thm:main}
Let $\Q$ be the field of rational numbers and $\Q_{\alg}$ a fixed algebraic closure of $\Q$. Let $F$ be the algebraic extension of $\Q$ obtained as the compositum in $\Q_{\alg}$ of all number fields in which each prime in $\{2,5,11,13\}$ splits completely. There is a cyclic division algebra $D$ of degree $6$ with center $F$ and a maximal subfield $K$ such that the following hold.
\begin{enumerate}[label=(\roman*)]
\item The normal closure of $K/F$ has Galois group $S_6$, and $N_{K/F}$ is surjective.
\item For every $x\in D^\times\setminus K^\times$, we have $\langle\ker N_{K/F},x(\ker N_{K/F})x^{-1}\rangle=D^{(1)}$.
\item The groups $K^\times$ and $\ker N_{K/F}$ are maximal subgroups of $D^\times$ and $D^{(1)}$, respectively.
\end{enumerate}
\end{theorem}
\end{samepage}

We construct the division algebra $D$ and its maximal subfield $K$ in \S\ref{sec:construction}, where the terminology on complete splitting is also recalled. By Proposition~\ref{prop:cyclic}, there are a cyclic maximal subfield $C$ of $D$ with $\Gal(C/F)=\langle\sigma\rangle$ and an element $a\in\Q^\times$ such that $D\simeq(C/F,\sigma,a)$. By contrast, the maximal subfield $K$ of Theorem~\ref{thm:main} is non-Galois over $F$, and $K^\times$ is an abelian maximal subgroup of $D^\times$.

Theorem~\ref{thm:main} should be compared with earlier results on maximal subgroups. For a noncommutative centrally finite division ring $D$, it follows from Akbari, Mahdavi-Hezavehi, and Mahmudi \cite[Theorem~7, p.~429]{AMM1999} that every nilpotent maximal subgroup of $D^\times$ is the multiplicative group of a maximal subfield. For almost subnormal subgroups of $D^\times$, it is shown in Khanh--Hai \cite[Theorems~1.1 and 1.2(i), p.~247]{KhanhHai2024} that every locally nilpotent maximal subgroup is abelian, and that a nonabelian locally solvable maximal subgroup can occur only if $D$ is cyclic of prime degree. The earlier result in Khanh \cite[Theorem~2.7(i)]{Khanh2020} assumes subnormality and algebraicity over the center of a term of the derived series of the maximal subgroup. In particular, every locally solvable maximal subgroup of $D^\times$ in our example is abelian, since $\deg D=6$.

In this example, a maximal subgroup $M$ of $D^\times$ is also abelian if it contains no noncyclic free subgroup or is solvable-by-finite. To see this, suppose $M$ is nonabelian. Then, it follows from Hai--Khanh \cite[Theorem~3.1, p.~610]{HaiKhanh2019} in the first case (with $D$ locally finite) and from Khanh--Hai \cite[Theorem~1.1, p.~78]{KhanhHai2022} in the second (with $|F|\geq5$) that $M$ has a simple quotient isomorphic to $\Gal(P/F)$ for a Galois maximal subfield $P$ of $D$. This is impossible, since $[P:F]=6$ and every group of order $6$ has a proper normal Sylow $3$-subgroup.

Recall that a subgroup $B$ of a group $H$ is \textit{malnormal} if $B\cap hBh^{-1}=1$ for every $h\in H\setminus B$. As shown in Khanh \cite[Proposition~5.2]{KhanhValued2026}, an abelian maximal subgroup of $D^\times$, for a noncommutative centrally finite division ring $D$ with center $F$, is $K^\times$ for a maximal subfield $K$, and $K^\times/F^\times$ is malnormal in $D^\times/F^\times$. Corollary~\ref{cor:structure} deduces malnormality for the pair $(D,K)$ of Theorem~\ref{thm:main} from the absence of proper intermediate fields and the equality $\Aut(K/F)=1$ of Proposition~\ref{prop:base-change}(ii). For these properties and the associated intersection criterion, see Hai--Huynh \cite[Corollary~4 and Theorem~1, pp.~40--41]{HaiHuynh2005}. According to Khanh \cite[Theorem~1.2]{KhanhValued2026}, the multiplicative group of a noncommutative centrally finite division ring has no abelian maximal subgroup if the center is Henselian for a nontrivial valuation and the ring is tame for the extended valuation. By contrast, every valuation on the division ring of Theorem~\ref{thm:main} is trivial; see Proposition~\ref{prop:valuations}.

The paper is organized as follows. In \S\ref{sec:construction}, we construct a division algebra $D_0$ with center $\Q$, with a maximal subfield $K_0$ and local indices $2$ and $3$, and we extend scalars from $\Q$ to $F$. Let $K_L=K_0\otimes_\Q L$ for any number field $L\subseteq F$. The surjectivity of $N_{K/F}$ is proved in \S\ref{sec:norm}: given $c\in F^\times$, we choose such an $L$ containing $c$. The equation $N_{K_L/L}(u)=c$ has a solution over each prescribed completion of $L$; hence, $c\in N_{K/F}(K^\times)$ by Lemma~\ref{lem:MB}.

The main step is carried out in \S\ref{sec:generation}. A reduced trace calculation shows that the multiplication map on suitable conjugates of the variety defined by $N_{K_L/L}=1$ has surjective differential at $(1,\ldots,1)$. Hence, over each prescribed completion, the group generated by these conjugates contains a neighborhood of $1$. For each elementary matrix $u$, a suitable diagonal matrix $a$ in this group satisfies $a^Nua^{-N}\to1$, from which the local factorizations of elements of reduced norm $1$ follow. Appending copies of $(1,\ldots,1)$ produces nonsingular local points of the corresponding fibers. By results on convolutions, these fibers are geometrically irreducible for sufficiently long products (see Lemma~\ref{lem:fibers}), so Lemma~\ref{lem:MB} applies. In \S\ref{sec:proof}, the resulting factorizations over number fields contained in $F$, together with the surjectivity of $N_{K/F}$, yield both maximality assertions of Theorem~\ref{thm:main}. Corollary~\ref{cor:structure} then uses the vanishing of the reduced Whitehead group and the generation assertion of Theorem~\ref{thm:main} to prove that $D^{(1)}/Z(D^{(1)})$ is nonabelian simple.

We use the following conventions throughout. For a field $L$, we write $L_{\alg}$ for an algebraic closure of $L$. The \textit{compositum} of a family of subfields of a common field is the subfield they generate; it is written $P\cdot L$ for two subfields $P$ and $L$. If $P/L$ is a field extension, we write $\Aut(P/L)$ for its group of $L$-automorphisms, and $\Gal(P/L)$ for this group when $P/L$ is Galois. If $P$ is a finite-dimensional commutative $L$-algebra, then the norm $N_{P/L}(u)$ and the trace $\tr_{P/L}(u)$ of an element $u\in P$ are the determinant and the trace of multiplication by $u$, respectively. We write $\Trd_A$ for the reduced trace of a central simple algebra $A$, and $\Br(L)$ for the Brauer group of $L$. Note that the reduced norm commutes with extension of scalars, by Draxl \cite[\S22, Theorem~2, p.~148]{Draxl1983}. The \textit{reduced Whitehead group} of a centrally finite division ring $D$ is $\SK_1(D)=D^{(1)}/[D^\times,D^\times]$ (see Tignol--Wadsworth \cite[Ch.~11, p.~535]{TignolWadsworth2015}).

\section{The division algebra \texorpdfstring{$D$}{D} and its maximal subfield \texorpdfstring{$K$}{K}}\label{sec:construction}

In this section, we construct the division algebra $D$ and the maximal subfield $K$ of Theorem~\ref{thm:main}. In \S\ref{ss:K0}, we choose a polynomial $f$ of degree $6$ with prescribed reductions modulo small primes; it defines a number field $K_0$, for which we determine the local structure at the primes $2$, $5$, $11$, $13$ and the Galois group of the normal closure. In \S\ref{ss:D0}, we prescribe local invariants to obtain a division algebra $D_0$ with center $\Q$ that contains $K_0$ as a maximal subfield. The field $F$ is studied in \S\ref{ss:F}. In \S\ref{ss:scalar}, we extend scalars from $\Q$ to the subfields of $F$; the resulting division algebra $D$ over $F$ is cyclic.

\subsection{The field \texorpdfstring{$K_0$}{K0}}\label{ss:K0}

We use the following terminology for valuations, as in Tignol--Wadsworth \cite[\S1.1.1, p.~2]{TignolWadsworth2015}. Let $\Gamma$ be a totally ordered abelian group, written additively. Adjoin to $\Gamma$ an element $\infty$ with $\gamma<\infty$ for $\gamma\in\Gamma$ and $\gamma+\infty=\infty+\gamma=\infty$ for $\gamma\in\Gamma\cup\{\infty\}$. A \textit{valuation} on a division ring $A$ is a map $v\colon A\to\Gamma\cup\{\infty\}$ satisfying the following conditions for all $x$, $y\in A$:
\begin{enumerate}[label=(\roman*)]
\item $v(x)=\infty$ if and only if $x=0$;
\item $v(xy)=v(x)+v(y)$;
\item $v(x+y)\ge\min\{v(x),v(y)\}$.
\end{enumerate}
Its \textit{valuation ring} is $\mathcal O_A=\{a\in A\mid v(a)\ge0\}$, with maximal ideal $\mathfrak m_A=\{a\in A\mid v(a)>0\}$. The quotient $\overline{A}=\mathcal O_A/\mathfrak m_A$ is the \textit{residue division ring} of $A$; when $A$ is a field, this quotient is its \textit{residue field}.

Let $L$ be a number field, and write $\mathbb Z_L$ for its ring of integers, the integral closure of $\mathbb Z$ in $L$. Let $p$ be a rational prime. We identify a \textit{place $v$ of $L$ above $p$} with a nonzero prime ideal $\mathfrak p_v$ of $\mathbb Z_L$ satisfying $\mathfrak p_v\cap\mathbb Z=p\mathbb Z$ (see Cohen \cite[\S4.6.2, p.~187, and Def.~4.8.2, p.~196]{Cohen1996}). For $x\in L^\times$, write $\operatorname{ord}_{\mathfrak p_v}(x)$ for the exponent of $\mathfrak p_v$ in the factorization of the fractional ideal $x\mathbb Z_L$ (see \cite[Theorem~4.6.14(3), p.~186]{Cohen1996}). The valuation $u_v$ on $L$ representing $v$ is defined by
$$
u_v(x)=\frac{\operatorname{ord}_{\mathfrak p_v}(x)}{\operatorname{ord}_{\mathfrak p_v}(p)}\quad\text{for }x\in L^\times\quad\text{and}\quad u_v(0)=\infty.
$$
If $a\in\mathbb Z$ and $p\nmid a$, then $a\notin\mathfrak p_v$, so $\operatorname{ord}_{\mathfrak p_v}(a)=0$. Thus $u_v|_\Q=v_p$, where $v_p$ is the $p$-adic valuation on $\Q$ normalized by $v_p(p)=1$. We write $L_v$ for the completion of $L$ with respect to $u_v$, constructed as in Jahnke \cite[\S3.1, pp.~134--135]{Jahnke2018}.

Following Lang \cite[Ch.~II, \S1, pp.~39--40]{Lang1994}, we say that $p$ \textit{splits completely} in $L$ if the canonical embedding $\Q_p\to L_v$ is an isomorphism for every place $v$ of $L$ above $p$. Equivalently, $p\mathbb Z_L$ is a product of $[L:\Q]$ distinct prime ideals \cite[Def.~4.8.7, p.~197]{Cohen1996}.

Let $P/k$ be a finite extension of $p$-adic fields. We write $\Gamma_k$ and $\overline{k}$ for the value group and the residue field of $k$, respectively, and we use the corresponding notation for $P$ with its extended valuation. Then, $[P:k]=ef$, where $e=|\Gamma_P:\Gamma_k|$ is the \textit{ramification index} and $f=[\overline{P}:\overline{k}]$ is the \textit{residue degree}. The extension is \textit{unramified} if $e=1$, \textit{ramified} if $e>1$, and \textit{totally ramified} if $f=1$, equivalently if $e=[P:k]$ (see Tignol--Wadsworth \cite[App.~A.1, Defs.~A.3--A.4, p.~586]{TignolWadsworth2015}).

The local data for our construction are chosen as follows. Throughout the paper, let
\[
\Sigma=\{2,5,11,13\}.
\]
We prescribe local indices $d_2=d_{11}=3$ and $d_5=d_{13}=2$. Their least common multiple is $6$, and Brauer invariants of these orders can be chosen in pairs summing to zero in $\Q/\mathbb Z$. Note that the local matrix sizes $m_p=6/d_p$ are at least $2$, and that $p\nmid d_p$ for $p\in\Sigma$.

At each $p\in\Sigma$, we seek $m_p$ field factors of degree $d_p$. Each such factor splits the local Brauer class, and the corresponding idempotents have rank one in the local matrix algebra. Moreover, the norm on units of an unramified factor maps onto $\mathbb Z_p^\times$, and the norm of a uniformizer in a totally ramified factor has $p$-adic valuation $1$. Therefore, for a monic polynomial $f\in\mathbb Z[X]$ of degree $6$, we prescribe
\begin{equation}\label{eq:prescribed}
f\bmod p=X^{d_p}\,\bar h_{p,2}\cdots\bar h_{p,m_p}\qquad\text{for }p\in\Sigma,
\end{equation}
with distinct monic irreducible polynomials $\bar h_{p,i}$ of degree $d_p$ and nonzero constant terms. With $f(0)=2\cdot5\cdot11\cdot13=1430$, the monic factor reducing to $X^{d_p}$ is then an Eisenstein polynomial, as we will see in the proof of Lemma~\ref{lem:local-K0}.

We also prescribe the reductions of $f$ modulo $3$ and $7$. The irreducible reduction at $7$ ensures transitivity of the Galois action. At $3$, the squarefree factorization of type $(1,5)$ makes Frobenius a $5$-cycle. At a place of the splitting field of $f$ above $5$, the ramified quadratic factor supplies a transposition in the inertia subgroup. We will use the same transposition in \S\ref{ss:F} to prove linear disjointness from $F$ over $\Q$. The reductions are as follows:
\begin{equation}\label{eq:table}
\begin{array}{c|l}
p & f(X)\bmod p\\ \hline
2  & X^3(X^3+X+1)\\
3  & (X-1)(X^5-X+1)\\
5  & X^2(X^2-2)((X-1)^2-2)\\
7  & X^6+2\\
11 & X^3(X^3+X+4)\\
13 & X^2(X^2-2)((X-1)^2-2)
\end{array}
\end{equation}
The constant terms agree with $1430$ modulo the respective primes. By the Chinese Remainder Theorem, the coefficients of $X,\ldots,X^5$ are determined modulo $30030=2\cdot3\cdot5\cdot7\cdot11\cdot13$. Taking their representatives in $[-15015,15015)$, we obtain
\begin{equation}\label{eq:f}
f(X)=X^6+14168X^5-7413X^4+12159X^3+5852X^2+10010X+1430.
\end{equation}

The factors displayed in \eqref{eq:table}, other than the powers of $X$, are irreducible. This is clear for the cubics, which have no roots in the respective residue fields, and for the quadratics, since $2$ is a nonsquare modulo $5$ and $13$. Consider now $h=X^5-X+1$ over $\Fq_3$. We have $X^9-X\equiv -X^4-1\pmod h$, hence
\[
\gcd(h,X^9-X)=\gcd(h,X^4+1)=\gcd(X+1,X^4+1)=1.
\]
But every irreducible polynomial over $\Fq_3$ of degree $1$ or $2$ divides $X^9-X=\prod_{u\in\Fq_9}(X-u)$, since its roots lie in $\Fq_9$. Therefore, each irreducible factor of $h$ has degree at least $3$, and $h$ must be irreducible because $\deg h=5<3+3$. It remains to consider $f\bmod 7=X^6+2$. Let $\alpha$ be a root of this polynomial. Since $\alpha^6=5$ has order $6$ in $\Fq_7^\times$, the order $n$ of $\alpha$ satisfies $n/\gcd(n,6)=6$, so $6\mid n$ and $n=36$. The conjugates of $\alpha$ over $\Fq_7$ are the distinct iterates $\alpha^{7^j}$ under $a\mapsto a^7$. Since $7^m\equiv1+6m\pmod{36}$, we have
\[
\alpha^{7^m}=\alpha\quad\Longleftrightarrow\quad 36\mid 7^m-1\quad\Longleftrightarrow\quad 6\mid m.
\]
Thus, $[\Fq_7(\alpha):\Fq_7]=\min\{m\ge1\mid\alpha^{7^m}=\alpha\}=6$. Since $f\bmod 7$ is monic of degree $6$, it is the minimal polynomial of $\alpha$ over $\Fq_7$. Therefore, $f\bmod7$ is irreducible, hence $f$ is irreducible over $\Q$.

We regard $K_0=\Q[X]/(f)$ as a subfield of $\Q_{\alg}$, and we let $E_0\subseteq\Q_{\alg}$ be its normal closure over $\Q$, which is the splitting field of $f$. The next lemma shows that $K_0$ has the prescribed local structure at the primes in $\Sigma$.

\begin{lemma}\label{lem:local-K0}
For $p\in\Sigma$, there is an isomorphism
\begin{equation}\label{eq:local-K0}
K_0\otimes_\Q\Q_p\simeq{\textstyle\prod_{i=1}^{m_p}}\,P_{p,i},
\end{equation}
where each $P_{p,i}$ is a field of degree $d_p$ over $\Q_p$. Exactly one factor is totally ramified, and the others are unramified.
\end{lemma}

\begin{proof}
Let $p\in\Sigma$, and let $\bar h_{p,2},\ldots,\bar h_{p,m_p}$ be the factors of $f\bmod p$ in \eqref{eq:table} other than $X^{d_p}$. These polynomials are distinct and irreducible, and $\bar h_{p,i}(0)\ne0$; hence, each factor of $f\bmod p$ in \eqref{eq:prescribed} is coprime to the product of the rest. Applying Hensel's lemma successively to such pairs over the complete ring $\mathbb Z_p$, we obtain monic polynomials $g_p,h_{p,2},\ldots,h_{p,m_p}\in\mathbb Z_p[X]$ with
\begin{equation}\label{eq:hensel}
f=g_ph_{p,2}\cdots h_{p,m_p},\quad\text{where}\quad
g_p\bmod p=X^{d_p}\quad\text{and}\quad
h_{p,i}\bmod p=\bar h_{p,i}\ \text{for }2\le i\le m_p.
\end{equation}
Since $g_p$ and each $h_{p,i}$ are monic, their degrees are the degrees $d_p$ of their reductions.

Consider first the factor $g_p$. Since $\bar h_{p,i}(0)\ne0$, we have $h_{p,i}(0)\in\mathbb Z_p^\times$ for $i\ge2$. In view of the reduction of $g_p$ and the factorization of $f$ in \eqref{eq:hensel}, it follows that
\[
g_p-X^{d_p}\in p\mathbb Z_p[X]\quad\text{and}\quad v_p(g_p(0))=v_p(f(0))=1.
\]
Thus, $g_p$ is an Eisenstein polynomial, hence $P_{p,1}=\Q_p[X]/(g_p)$ is a totally ramified extension of $\Q_p$ of degree $d_p$.

Now, let $i\ge2$. Since $\bar h_{p,i}$ is irreducible, $h_{p,i}$ is irreducible over $\Q_p$; hence, $P_{p,i}=\Q_p[X]/(h_{p,i})$ is a field of degree $d_p$ over $\Q_p$. The class $\beta_i$ of $X$ in $P_{p,i}$ is integral over $\mathbb Z_p$ since $h_{p,i}$ is monic, and $\bar h_{p,i}(\bar\beta_i)=0$. Therefore,
\[
d_p=[\Fq_p(\bar\beta_i):\Fq_p]\le[\overline{P_{p,i}}:\Fq_p]\le[P_{p,i}:\Q_p]=d_p.
\]
Thus, equality holds throughout, and $P_{p,i}/\Q_p$ is unramified. Since the factors in \eqref{eq:hensel} are pairwise coprime, the Chinese Remainder Theorem yields the isomorphism \eqref{eq:local-K0}.
\end{proof}

To determine the Galois group of $E_0$ over $\Q$ in Lemma~\ref{lem:galois} below, we use decomposition and inertia groups, for which we recall the definitions from Tignol--Wadsworth \cite[App.~A.1, p.~589]{TignolWadsworth2015}. Let $M/L$ be a Galois extension of number fields, and let $w$ be a place of $M$ above a rational prime $p$. Let $\mathrm{P}_w$ be the corresponding prime ideal of $\mathbb Z_M$ and $v$ the place of $L$ corresponding to $\mathrm{P}_w\cap\mathbb Z_L$ (see \cite[Ch.~II, \S1, pp.~37--38]{Lang1994}). If $u$ is the valuation on $M$ representing $w$ with $u(p)=1$, then $u|_L=u_v$. The \textit{decomposition group at $w$} is
$$
G_w=\{\sigma\in\Gal(M/L)\mid u\circ\sigma=u\}=\{\sigma\in\Gal(M/L)\mid\sigma(\mathrm{P}_w)=\mathrm{P}_w\}.
$$
Each $\sigma\in G_w$ extends to an $L_v$-automorphism of $M_w$ and induces an automorphism $\bar\sigma$ of $\overline{M_w}$ by $\bar\sigma(\bar a)=\overline{\sigma(a)}$ for $a\in\mathcal O_{M_w}$. The \textit{inertia subgroup at $w$} is the kernel of the reduction map $\rho_w\colon G_w\to\Gal(\overline{M_w}/\overline{L_v})$ defined by $\sigma\mapsto\bar\sigma$. Thus,
\[
I_w=\ker\rho_w=\{\sigma\in G_w\mid \sigma(a)-a\in\mathfrak m_{M_w}\text{ for all }a\in\mathcal O_{M_w}\}.
\]

\begin{lemma}\label{lem:galois}
The action on the roots of $f$ identifies $\Gal(E_0/\Q)$ with $S_6$, and the inertia subgroup at a place of $E_0$ above $5$ contains a transposition.
\end{lemma}

\begin{proof}
Let $G=\Gal(E_0/\Q)$, and let $\Omega$ be the set of roots of $f$ in $E_0$. Since $f$ is irreducible over $\Q$, the group $G$ acts transitively on $\Omega$. We will show that $G$ contains a $5$-cycle and that the inertia subgroup at a place of $E_0$ above $5$ contains a transposition, and then derive the equality $G=S_6$.

Take any place $w_3$ of $E_0$ above $3$, and consider reductions in $\overline{(E_0)_{w_3}}$. Let $h=X^5-X+1\in\Fq_3[X]$, as above. The roots of $f$ are integral at $w_3$ since $f$ is monic, and $\prod_{u\in\Omega}(X-\bar u)=(X-1)h$ is squarefree. Thus, reduction identifies the two root sets. Since every element of $I_{w_3}$ fixes each $\bar u$, it fixes each $u\in\Omega$; hence, $I_{w_3}=1$. The Frobenius element of $G_{w_3}$ acts on the reduced roots by $\bar u\mapsto\bar u^3$, with orbit lengths $1$ and $5$, the degrees of $X-1$ and $h$. Therefore, it acts on $\Omega$ as a $5$-cycle.

Now, let $w_5$ be a place of $E_0$ above $5$, and let $U$ be the maximal unramified subextension of $(E_0)_{w_5}/\Q_5$. Let $g_5$, $h_{5,2}$, $h_{5,3}$ be the factors of $f$ over $\mathbb Z_5$ given by \eqref{eq:hensel}. If $u$ is a root of $h_{5,2}$ or $h_{5,3}$, then the extension $\Q_5(u)/\Q_5$ is unramified, as shown in the proof of Lemma~\ref{lem:local-K0}; hence, $u\in U$. Since $U/\Q_5$ is unramified, $5$ remains a uniformizer of $U$, and
\[
g_5\bmod 5=X^2\quad\text{in }\overline U[X]\quad\text{and}\quad g_5(0)\in5\mathcal O_U\setminus5^2\mathcal O_U.
\]
Thus, $g_5$ is an Eisenstein polynomial over $U$. Its roots $b$ and $c$ satisfy $c\in U(b)$, since $g_5$ is quadratic. As $(E_0)_{w_5}$ is the splitting field of $f$ over $\Q_5$, it follows that
\[
(E_0)_{w_5}=U(b)\quad\text{and}\quad[U(b):U]=2.
\]
Therefore, the nonidentity element of $I_{w_5}\simeq\Gal((E_0)_{w_5}/U)$ fixes the roots of the $h_{5,i}$ and interchanges $b$ and $c$, so it acts on $\Omega$ as the transposition $\tau=(b\,c)$.

To complete the proof, let $a$ be the fixed point of the $5$-cycle. Its stabilizer $G_a$ is transitive on $\Omega\setminus\{a\}$. Since $G$ is transitive on $\Omega$, it must be doubly transitive, hence
\[
\{g\tau g^{-1}\mid g\in G\}=\{(g(b)\,g(c))\mid g\in G\}=\{(u\,v)\mid u,v\in\Omega,\ u\ne v\}.
\]
These transpositions generate $S_6$, so $G=S_6$.
\end{proof}

\subsection{The division algebra \texorpdfstring{$D_0$}{D0}}\label{ss:D0}

In this subsection, we define a division algebra with center $\Q$ by prescribing its local invariants. We use the following terminology. Let $A$ be a central simple algebra over a field $L$. Its dimension is denoted by $[A:L]$, and its degree is $\deg A=\sqrt{[A:L]}$. Write $A\simeq M_r(\Delta)$, where $\Delta$ is the division algebra associated to $A$, and let $\ind A=\deg\Delta$. We use the same index notation for Brauer classes. If $k$ is a completion of a number field, then $\inv_k\colon\Br(k)\to\Q/\mathbb Z$ denotes the local invariant map, defined in Tignol--Wadsworth \cite[\S1.2.8, pp.~23--24]{TignolWadsworth2015} for nonarchimedean $k$. At a real place, the nonzero class has invariant $1/2+\mathbb Z$, and at a complex place the Brauer group is trivial. We write $\beta_L$ for the restriction of a class $\beta\in\Br(\Q)$ to an extension $L/\Q$, and we let $\inv_v(\beta)=\inv_{\Q_v}(\beta_{\Q_v})$ for a place $v$ of $\Q$. By the restriction formula of Serre \cite[Chap.~XIII, \S3, Prop.~7, pp.~193--194]{Serre}, for a rational prime $p$ and a finite extension $P/\Q_p$ we have
\begin{equation}\label{eq:restriction}
\inv_{P}(\beta_{P})=[P:\Q_p]\inv_p(\beta).
\end{equation}

Since the invariants prescribed in \S\ref{ss:K0} sum to zero in $\Q/\mathbb Z$, the Albert--Brauer--Hasse--Noether Theorem (see H\"urlimann \cite[Theorem~4, p.~219]{Hurlimann}) ensures that there is a class $\beta\in\Br(\Q)$ with $\inv_v(\beta)=0$ for $v\notin\Sigma$ and
\begin{equation}\label{eq:invariants}
(\inv_2(\beta),\inv_5(\beta),\inv_{11}(\beta),\inv_{13}(\beta))=(1/3,1/2,-1/3,1/2).
\end{equation}
Let $D_0$ be the division algebra in the class $\beta$. Its basic properties are as follows:

\begin{proposition}\label{prop:base-division}
The algebra $D_0$ is a centrally finite division ring of degree $6$ with center $\Q$. Moreover, the following hold.
\begin{enumerate}[label=(\roman*)]
\item For $p\in\Sigma$, we have
\begin{equation}\label{eq:local-D0}
D_0\otimes_\Q\Q_p\simeq M_{m_p}(\Delta_p),
\end{equation}
where $\Delta_p$ is a division ring of degree $d_p$ with center $\Q_p$.
\item The algebra $D_0$ splits at every place of $\Q$ outside $\Sigma$.
\item The field $K_0$ embeds in $D_0$, and its image under any embedding is a maximal subfield of $D_0$.
\end{enumerate}
\end{proposition}

\begin{proof}
For $p\in\Sigma$, the factors $P_{p,i}$ of \eqref{eq:local-K0} satisfy
\[
\inv_{P_{p,i}}(\beta_{P_{p,i}})=[P_{p,i}:\Q_p]\inv_p(\beta)=d_p\inv_p(\beta)=0
\]
by \eqref{eq:restriction} and \eqref{eq:invariants}. Moreover, the invariants of $\beta_{K_0}$ at the other places of $K_0$ vanish since $\inv_v(\beta)=0$ for $v\notin\Sigma$. Hence, $\beta_{K_0}=0$ by the Albert--Brauer--Hasse--Noether Theorem, i.e., $K_0$ splits $\beta$.

For $p\in\Sigma$, let $\Delta_p$ be the division algebra in the class $\beta_{\Q_p}$, with $n_p=\deg\Delta_p=\ind(\beta_{\Q_p})$, and let $M_p$ be a maximal subfield of $\Delta_p$. Since $P_{p,i}$ and $M_p$ split $\Delta_p$, index divisibility and the restriction formula \eqref{eq:restriction} imply
\[
n_p\mid[P_{p,i}:\Q_p]=d_p\quad\text{and}\quad
0=\inv_{M_p}(\beta_{M_p})=[M_p:\Q_p]\inv_p(\beta)=n_p\inv_p(\beta).
\]
But $\inv_p(\beta)$ has order $d_p$, so $d_p\mid n_p$, hence $n_p=d_p$. By scalar extension, each $n_p$ divides $\ind(\beta)$, and we have just seen that $K_0$ splits $\beta$. Hence,
\[
6=\operatorname{lcm}_{p\in\Sigma} n_p\mid\ind(\beta)\mid[K_0:\Q]=6.
\]
This shows that $\deg D_0=\ind(\beta)=6$. For $p\in\Sigma$, the matrix size of $D_0\otimes_\Q\Q_p$ is therefore $\deg D_0/\deg\Delta_p=6/d_p=m_p$, proving (i). If $v$ is a place of $\Q$ outside $\Sigma$, then $D_0$ splits at $v$ since $\inv_v(\beta)=0$, proving (ii).

For (iii), note that via $D_0\hookrightarrow D_0\otimes_\Q K_0\simeq M_6(K_0)$, the space $V=K_0^6$ is a left $D_0$-module of rank one, since $\dim_\Q V=6[K_0:\Q]=36=[D_0:\Q]$. Scalar multiplication commutes with the matrix action, so we have the embedding
\[
K_0\hookrightarrow\End_{D_0}V\quad\text{defined by}\quad a\longmapsto (v\mapsto av).
\]
Right multiplication on the regular left module $D_0\simeq V$ identifies $D_0^{\mathrm{op}}$ with $\End_{D_0}V$. Taking opposite rings embeds $K_0=K_0^{\mathrm{op}}$ into $D_0$. The image of any such embedding is a maximal subfield since $[K_0:\Q]=\deg D_0$.
\end{proof}

\subsection{The field \texorpdfstring{$F$}{F}}\label{ss:F}

Let $F$ be the field of Theorem~\ref{thm:main}, i.e., the compositum in $\Q_{\alg}$ of all number fields in which each prime in $\Sigma$ splits completely. Complete splitting at the primes in $\Sigma$ will preserve the local indices of Proposition~\ref{prop:base-division} over any number field $L\subseteq F$. At $5$, it also forces inertia in $\Gal(E_0/\Q)$ to fix $F\cap E_0$. We will need the following properties of $F$:

\begin{lemma}\label{lem:center}
\begin{enumerate}[label=(\roman*)]
\item Every prime in $\Sigma$ splits completely in every number field $L\subseteq F$.
\item The extension $F/\Q$ is Galois, and $F\cap E_0=\Q$.
\end{enumerate}
\end{lemma}

\begin{proof}
(i) By the criterion in Lang \cite[Ch.~II, \S1, pp.~39--40]{Lang1994}, a prime $p$ splits completely in a number field $L$ if and only if $\sigma(L)\subseteq\Q_p$ for every embedding $\sigma\colon L\hookrightarrow(\Q_p)_{\alg}$. Let $L\subseteq F$ be a number field, and choose fields $L_1,\ldots,L_r$ from the definition of $F$ with $L\subseteq L'=L_1\cdot\ldots\cdot L_r$. Let $p\in\Sigma$, and let $\sigma\colon L\hookrightarrow(\Q_p)_{\alg}$ be an embedding. Extend $\sigma$ to $\widetilde\sigma\colon L'\hookrightarrow(\Q_p)_{\alg}$. Then,
\[
\sigma(L)\subseteq\widetilde\sigma(L')=\widetilde\sigma(L_1)\cdot\ldots\cdot\widetilde\sigma(L_r)\subseteq\Q_p.
\]
Therefore, every $p\in\Sigma$ splits completely in $L$ (see Lang \cite[Ch.~II, \S1, SC~1 and SC~3, p.~40]{Lang1994}).

(ii) The criterion in (i) is invariant under $\Q$-conjugation, so the fields defining $F$ are permuted by $\Aut(\Q_{\alg}/\Q)$. Hence, their compositum $F$ is Galois over $\Q$.

Now, let $M=F\cap E_0$. By Lemma~\ref{lem:galois}, we may choose a place $w$ of $E_0$ above $5$ and a transposition $\tau\in I_w$. Let $\mathrm{P}_w$ be the prime ideal of $\mathbb Z_{E_0}$ corresponding to $w$, and let $v$ be the place of $M$ corresponding to $\mathrm{P}_w\cap\mathbb Z_M$. Since $5$ splits completely in $M$ by (i), we have
\[
M\subseteq M_v=\Q_5\subseteq(E_0)_w.
\]
Each $\sigma\in I_w$ acts on $(E_0)_w$ over $\Q_5$, hence $\sigma|_M=\mathrm{id}_M$. Moreover, $M/\Q$ is Galois since $F/\Q$ and $E_0/\Q$ are Galois. Thus,
\[
\tau\in I_w\subseteq\Gal(E_0/M)\trianglelefteq\Gal(E_0/\Q)=S_6.
\]
But the $S_6$-conjugates of $\tau$ are all the transpositions. Hence, by normality,
\[
S_6=\langle g\tau g^{-1}\mid g\in S_6\rangle\subseteq\Gal(E_0/M)\subseteq S_6.
\]
Therefore, $M=E_0^{S_6}=\Q$.
\end{proof}

\subsection{Scalar extension to subfields of \texorpdfstring{$F$}{F}}\label{ss:scalar}

We fix an embedding $K_0\hookrightarrow D_0$ as in Proposition~\ref{prop:base-division}(iii), and we view $K_0$ as a maximal subfield of $D_0$. If $L$ is a field with $\Q\subseteq L\subseteq F$, let
\begin{equation}\label{eq:DL}
D_L=D_0\otimes_\Q L\quad\text{and}\quad K_L=K_0\otimes_\Q L\subseteq D_L.
\end{equation}
We collect in the following proposition the properties of $D_L$ and $K_L$ that are used in the rest of the paper.

\begin{proposition}\label{prop:base-change}
Let $L\subseteq F$ be a field with $L=F$ or $[L:\Q]<\infty$.
\begin{enumerate}[label=(\roman*)]
\item The algebra $D_L$ is a division ring of degree $6$ with center $L$, and $K_L$ is a maximal subfield of $D_L$.
\item The normal closure of $K_L/L$ is $E_0\cdot L$, and its Galois group acts as $S_6$ on the six embeddings of $K_L$ over $L$. The extension $K_L/L$ has no proper intermediate fields, and $\Aut(K_L/L)=1$.
\item Suppose $[L:\Q]<\infty$, and let $v$ be a place of $L$ above a prime $p\in\Sigma$. Then, the canonical embedding $\Q_p\to L_v$ is an isomorphism, and it yields identifications
\begin{equation}\label{eq:local-L}
D_L\otimes_L L_v=D_0\otimes_\Q\Q_p\simeq M_{m_p}(\Delta_p)\quad\text{and}\quad
K_L\otimes_L L_v=K_0\otimes_\Q\Q_p\simeq{\textstyle\prod_{\ell=1}^{m_p}}\,P_{p,\ell},
\end{equation}
where the isomorphisms are those of \eqref{eq:local-D0} and \eqref{eq:local-K0}.
\end{enumerate}
\end{proposition}

\begin{proof}
We first prove (iii), which is used in the proof of (i). By Lemma~\ref{lem:center}(i), the canonical embedding $\Q_p\to L_v$ is an isomorphism; hence, $D_L\otimes_L L_v=D_0\otimes_\Q L_v=D_0\otimes_\Q\Q_p$, and likewise for $K_L$. The isomorphisms in \eqref{eq:local-L} are then given by Proposition~\ref{prop:base-division}(i) and Lemma~\ref{lem:local-K0}.

(i) Clearly, $D_L$ is a central simple algebra of degree $6$ over $L$. Suppose first that $[L:\Q]<\infty$, and let $\Delta$ be the division algebra associated to $D_L$. If $v$ is a place of $L$ above $p\in\Sigma$, then \eqref{eq:local-L} shows that $\Delta\otimes_L L_v\simeq M_{r_v}(\Delta_p)$ for some $r_v\ge1$, since scalar extension preserves Brauer equivalence. Thus, $\ind D_L=\deg\Delta=r_vd_p$, so $d_p\mid\ind D_L$ for each $p\in\Sigma$, and
\[
6=\operatorname{lcm}_{p\in\Sigma}d_p\mid\ind D_L\mid\deg D_L=6.
\]
Therefore, the matrix size $\deg D_L/\ind D_L$ is $1$, i.e., $D_L$ is a division ring. Suppose now that $L=F$. Each nonzero $a\in D_F$ lies in $D_{L'}$ for some number field $L'\subseteq F$, hence $a^{-1}\in D_{L'}\subseteq D_F$. So, $D_F$ is a division ring. The assertion on $K_L$ will be proved together with (ii).

(ii) Let $E=E_0\cdot L$, which is the splitting field of $f$ over $L$. The restriction map $\rho\colon\Gal(E/L)\to\Gal(E_0/\Q)$ defined by $\sigma\mapsto\sigma|_{E_0}$ is injective since $E$ is generated by $E_0$ and $L$. By Lemma~\ref{lem:center}(ii), the fixed field of its image is
\[
E_0^{\operatorname{im}\rho}=E_0\cap E^{\Gal(E/L)}=E_0\cap L=\Q.
\]
Hence, by the Galois correspondence and Lemma~\ref{lem:galois},
\[
\Gal(E/L)\xrightarrow[\sim]{\ \rho\ }\Gal(E_0/\Q)=S_6\quad\text{and}\quad [E:L]=[E_0:\Q].
\]
The degree equality shows that $E_0$ and $L$ are linearly disjoint over $\Q$. Thus, $K_L\simeq K_0\cdot L$ is a field, with $[K_L:L]=[K_0:\Q]=6=\deg D_L$, hence it is a maximal subfield of $D_L$. This completes the proof of (i). The normal closure of $K_L/L$ is $E$, and the six $L$-embeddings of $K_L$ correspond to the roots of $f$, on which $\rho$ preserves the action. Since the point stabilizer $S_5$ is maximal in $S_6$, the extension $K_L/L$ has no proper intermediate fields. Also, the normalizer of $S_5$ fixes the unique point fixed by $S_5$, hence $N_{S_6}(S_5)=S_5$ and $\Aut(K_L/L)\simeq N_{S_6}(S_5)/S_5=1$.
\end{proof}

Henceforth, we let
\[
D=D_F\quad\text{and}\quad K=K_F.
\]
By Proposition~\ref{prop:base-change}, $D$ is a division ring of degree $6$ with center $F$, and $K$ is a maximal subfield of $D$. The next proposition shows that $D$ is a cyclic algebra; we use the notation of Lam \cite[\S14, pp.~218--219]{Lam2001} for cyclic algebras.

\begin{proposition}\label{prop:cyclic}
There are a cyclic maximal subfield $C$ of $D$, a generator $\sigma$ of $\Gal(C/F)$, and an element $a\in\Q^\times$ such that
\begin{equation}\label{eq:cyclic}
D\simeq(C/F,\sigma,a)={\textstyle\bigoplus_{i=0}^{5}}\,Cu^i,\quad\text{with}\quad u^6=a\quad\text{and}\quad uc=\sigma(c)u\ \text{for }c\in C.
\end{equation}
\end{proposition}

\begin{proof}
By Albert--Hasse \cite[p.~722]{AlbertHasse1932}, we may write $D_0\simeq(C_0/\Q,\sigma_0,a)$, where $C_0/\Q$ is cyclic of degree $6$, $\sigma_0$ generates $\Gal(C_0/\Q)$, and $a\in\Q^\times$. Let $C=C_0\otimes_\Q F$. Since $C$ embeds in $D$, it is a commutative domain of dimension $6$ over $F$, hence a field. Moreover, $C$ is a maximal subfield of $D$ since $[C:F]=\deg D$. The automorphism $\sigma=\sigma_0\otimes 1$ of $C$ has order $6$; hence, $\Gal(C/F)=\langle\sigma\rangle$. Therefore, extending scalars from $\Q$ to $F$ yields the isomorphism \eqref{eq:cyclic}.
\end{proof}

By Proposition~\ref{prop:base-change}(ii), the maximal subfield $K$ is non-Galois over $F$, and its normal closure has Galois group $S_6$. Theorem~\ref{thm:main}(iii) asserts that the abelian group $K^\times$ is maximal in $D^\times$.

\section{The norm from \texorpdfstring{$K$}{K} to \texorpdfstring{$F$}{F}}\label{sec:norm}

In this section, we prove that the norm $N_{K/F}\colon K^\times\to F^\times$ is surjective; see Proposition~\ref{prop:norm}. The proof is based on a local--global principle for varieties over the number fields contained in $F$, which we establish in \S\ref{ss:MB}. This principle is used again in \S\ref{sec:generation}.

\subsection{A local--global principle}\label{ss:MB}

We say that a place of a number field $L$ is \textit{above $\Sigma$} if it is above some prime $p\in\Sigma$. If $X$ is an affine variety over a field $L$ and $M/L$ is a field extension, we write $X(M)$ for the set of points of $X$ with coordinates in $M$, and $X_M$ for the scalar extension of $X$ to $M$. We use the following lemma to pass from local solutions at the places above $\Sigma$ to solutions over number fields contained in $F$.

\begin{lemma}\label{lem:MB}
Let $L\subseteq F$ be a number field, and let $X$ be an affine variety over $L$ that remains irreducible over $L_{\alg}$. If $X$ has a nonsingular point over $L_v$ for every place $v$ of $L$ above $\Sigma$, then $X$ has a point over some finite extension $L'/L$ contained in $F$.
\end{lemma}

\begin{proof}
Let $U=X_{\mathrm{sm}}$ be the smooth locus of $X$ over $L$. We first show that $U$ is geometrically connected and has a point over $L_v$ for every place $v$ of $L$ above $\Sigma$. Take any such $v$. Let $x_v\in X_{L_v}(L_v)$ be a nonsingular point, and let $A_v=\mathcal O_{X_{L_v},x_v}$, with maximal ideal $\mathfrak m_v$ and residue field $L_v$. Nonsingularity means that $A_v$ is regular, i.e.,
\[
\dim_{L_v}(\mathfrak m_v/\mathfrak m_v^2)=\dim A_v.
\]
Since $L_v$ has characteristic zero, it is perfect, and regularity at $x_v$ is equivalent to smoothness over $L_v$. If $M/L$ is an extension with projection $\pi_M\colon X_M\to X$, then $(X_M)_{\mathrm{sm}}=\pi_M^{-1}(U)=U_M$, because smoothness is preserved and detected by extension of the ground field. Taking $M=L_v$, we obtain
\[
x_v\in (X_{L_v})_{\mathrm{sm}}(L_v)=U_{L_v}(L_v)=U(L_v).
\]
Furthermore, the nonempty open subvariety $U_{L_{\alg}}\subseteq X_{L_{\alg}}$ is irreducible, hence $U$ is smooth and geometrically connected.

We may therefore apply the theorem of Moret-Bailly \cite[Theorem~1.3, p.~182]{MB} in the form of Barnet-Lamb et al.\ \cite[Proposition~3.1.1, pp.~545--546]{BLGGT}. In the notation of that proposition, we take $T=U$ and $K_0=K=K^{(\mathrm{avoid})}=L$, and we let $S$ be the set of places of $L$ above $\Sigma$. For $v\in S$, we prescribe the extension $L_v/L_v$ and the nonempty open set $\Omega_v=U(L_v)$. All prescribed extensions are trivial, so the Galois, compatibility, and invariance conditions hold. The resulting finite extension $L'/L$ satisfies
\[
U(L')\ne\varnothing\quad\text{and}\quad[L'_w:L_v]=1\quad\text{for every }w\mid v\text{ with }v\in S.
\]
Embed $L'$ into $\Q_{\alg}$ over $L$. If $p\in\Sigma$ and $w\mid v\mid p$, then $[L_v:\Q_p]=1$ by Lemma~\ref{lem:center}(i), so
\[
[L'_w:\Q_p]=[L'_w:L_v][L_v:\Q_p]=1.
\]
Therefore, every $p\in\Sigma$ splits completely in $L'$, hence $L'\subseteq F$ by the definition of $F$. Moreover, $X(L')\supseteq U(L')\ne\varnothing$.
\end{proof}

\subsection{Surjectivity of the norm}\label{ss:norm}

Locally at the primes in $\Sigma$, the surjectivity of the norm follows from Lemma~\ref{lem:local-K0}. Recall that, for $p\in\Sigma$, the algebra $K_0\otimes_\Q\Q_p$ is a product of fields of degree $d_p$ over $\Q_p$, one of which is totally ramified while the others are unramified. Its norm is the product of the norms of these fields.

\begin{lemma}\label{lem:local-norm}
For $p\in\Sigma$, the norm $N_{(K_0\otimes_\Q\Q_p)/\Q_p}$ maps $(K_0\otimes_\Q\Q_p)^\times$ onto $\Q_p^\times$.
\end{lemma}

\begin{proof}
By Lemma~\ref{lem:local-K0}, the algebra $K_0\otimes_\Q\Q_p$ has an unramified factor $U$ and a totally ramified factor $R$, both of degree $d_p$ over $\Q_p$. The norm $N_{U/\Q_p}$ maps $\mathcal O_U^\times$ onto $\mathbb Z_p^\times$. To see this, let $u\in\mathbb Z_p^\times$. The norm on $\mathcal O_U^\times$ reduces to the norm on residue field units, which is surjective, so there is some $z\in\mathcal O_U^\times$ with $u/N_{U/\Q_p}(z)\in1+p\mathbb Z_p$. Since $p\nmid d_p$, Hensel's lemma yields $t\in1+p\mathbb Z_p$ with $t^{d_p}=u/N_{U/\Q_p}(z)$. Then,
\[
N_{U/\Q_p}(zt)=N_{U/\Q_p}(z)t^{d_p}=u.
\]

Consider now a root $\alpha$ in $R$ of the Eisenstein polynomial $g_p$ of \eqref{eq:hensel}, and let $\eta=N_{R/\Q_p}(\alpha)=(-1)^{d_p}g_p(0)$. We have $v_p(\eta)=1$ by the Eisenstein condition. With $1$ in the other factors, the norm of $K_0\otimes_\Q\Q_p$ therefore takes all values in
\[
N_{U/\Q_p}(\mathcal O_U^\times)N_{R/\Q_p}(\langle\alpha\rangle)
=\mathbb Z_p^\times\eta^{\mathbb Z}
=\Q_p^\times.\qedhere
\]
\end{proof}

We can now prove the main result of this section.

\begin{proposition}\label{prop:norm}
The norm $N_{K/F}\colon K^\times\to F^\times$ is surjective.
\end{proposition}

\begin{proof}
Let $c\in F^\times$, and let $L=\Q(c)$, which is a number field contained in $F$. Let $X_c$ be the variety defined by the equation $N_{K_L/L}(u)=c$ in the affine space underlying $K_L$. We verify the hypotheses of Lemma~\ref{lem:MB} for $X_c$.

Over $L_{\alg}$, the six embeddings of $K_L$ identify the norm equation with $u_1\cdots u_6=c$. Projection onto the first five coordinates then identifies $(X_c)_{L_{\alg}}$ with
\[
\{(u_1,\ldots,u_5)\in\mathbb A^5_{L_{\alg}}\mid u_1\cdots u_5\ne0\},
\quad\text{with}\quad u_6=c/(u_1\cdots u_5).
\]
Thus, $X_c$ remains irreducible over $L_{\alg}$, and it is nonsingular. If $v$ is a place of $L$ above $p\in\Sigma$, then by Proposition~\ref{prop:base-change}(iii) we may identify $L_v$ with $\Q_p$ and $K_L\otimes_L L_v$ with $K_0\otimes_\Q\Q_p$. Hence, $X_c$ has a point over $L_v$ by Lemma~\ref{lem:local-norm}, and this point is nonsingular.

Therefore, by Lemma~\ref{lem:MB} the equation $N_{K_L/L}(u)=c$ has a solution $a\in K_{L'}^\times$ for some finite extension $L'/L$ inside $F$. As the norm is the determinant of multiplication, it commutes with scalar extension; hence, $N_{K/F}(a\otimes1)=N_{K_{L'}/L'}(a)=c$.
\end{proof}

\begin{remark}\label{rem:seven}
For comparison with Proposition~\ref{prop:norm}, consider the norm $N_{K_0/\Q}\colon K_0^\times\to\Q^\times$. The reduction $f\bmod 7$ is irreducible by \S\ref{ss:K0}, so $K_0\otimes_\Q\Q_7$ is an unramified field extension of degree $6$. Hence $\mathfrak p_7=7\mathbb Z_{K_0}$ is a prime ideal of residue degree $6$. By \cite[Prop.~4.6.15, p.~187, and Theorem~4.8.5, p.~197]{Cohen1996}, we have
$$
v_7\bigl(N_{K_0/\Q}(a)\bigr)=6\operatorname{ord}_{\mathfrak p_7}(a)\in 6\mathbb Z\quad\text{for }a\in K_0^\times.
$$
Since $v_7(7)=1$, the image of $N_{K_0/\Q}$ omits $7$.
\end{remark}

\section{Generation of the reduced norm-one group}\label{sec:generation}

Throughout this section, we fix a number field $L\subseteq F$ and an element $x\in D_L^\times\setminus K_L^\times$. Let $T$ and $V$ be the affine varieties defined by the equations $N_{K_L/L}=1$ and $\Nrd_{D_L}=1$ in the affine spaces underlying $K_L$ and $D_L$, respectively. For each field extension $M/L$, their points are described by
\begin{equation}\label{eq:TV}
T(M)=\{t\in K_L\otimes_L M\mid N_{(K_L\otimes_L M)/M}(t)=1\}
\quad\text{and}\quad V(M)=(D_L\otimes_L M)^{(1)}.
\end{equation}
Since $K_L$ is a maximal subfield of $D_L$ by Proposition~\ref{prop:base-change}(i), the reduced norm on $D_L$ restricts to the field norm on $K_L$ (see Tignol--Wadsworth \cite[proof of Theorem~1.4, p.~12]{TignolWadsworth2015}). By extension of scalars, we have $T(M)\subseteq V(M)$ for every $M/L$. Let
\begin{equation}\label{eq:H}
H=\langle T(L),xT(L)x^{-1}\rangle\subseteq D_L^{(1)}.
\end{equation}
The main result of this section is the following proposition, which is the key step in the proof of Theorem~\ref{thm:main}(ii):

\begin{proposition}\label{prop:factor}
For every $s\in D_L^{(1)}$, there is a finite extension $L'/L$ inside $F$ such that $s\in\langle T(L'),xT(L')x^{-1}\rangle$ in $D_{L'}^{(1)}$.
\end{proposition}

The proof is given in \S\ref{ss:proof-factor}, after some preparation. We consider products of elements taken from finitely many conjugates $h_jTh_j^{-1}$ of $T$, with $h_j\in H$. In \S\ref{ss:diff}, we choose the $h_j$ so that the corresponding multiplication map $\phi$ has surjective differential at $(1,\ldots,1)$. It is shown in \S\ref{ss:local} that the groups $h_jT(L_v)h_j^{-1}$ generate $(D_L\otimes_L L_v)^{(1)}$ for each place $v$ of $L$ above $\Sigma$, which yields local factorizations of $s$. Factoring $s$ over an extension of $L$ amounts to finding a point on a fiber of a convolution power of $\phi$. These fibers are geometrically irreducible by the results of \S\ref{ss:fibers}, so that Lemma~\ref{lem:MB} applies.

\subsection{The differential of the multiplication map}\label{ss:diff}

Let
\begin{equation}\label{eq:W}
K_L^0=\ker\tr_{K_L/L}\quad\text{and}\quad
W=\sum_{g\in H}gK_L^0g^{-1}\subseteq D_L.
\end{equation}
Note that $W$ is the $L$-subspace of $D_L$ spanned by the conjugates of $K_L^0$ under the elements of $H$. The computation of the differential in Lemma~\ref{lem:dense} below relies on the following description of $W$:

\begin{lemma}\label{lem:W}
Let $W$ be as in \eqref{eq:W}. Then, $W=\ker\Trd_{D_L}$.
\end{lemma}

\begin{proof}
Let $k=L_{\alg}$, and identify $A=D_L\otimes_L k$ with $M_6(k)$ so that $K_L\otimes_L k$ is the algebra of diagonal matrices. For $u\in K_L$, write $u_1,\ldots,u_6$ for its diagonal entries. Taking matrix traces in $A$, we see that $W\subseteq\ker\Trd_{D_L}$. To prove the reverse inclusion, it suffices to show that $W_k=W\otimes_L k$ contains every matrix of trace $0$ in $A$.

The first step is to choose a suitable element $t\in T(L)$. Let $u\in K_L^\times$ be such that the elements $(u_i/u_j)^6$ with $i\ne j$ are distinct and different from $1$. Such a $u$ exists: after clearing denominators, these conditions amount to avoiding the zeros of finitely many nonzero polynomials in the independent coordinates $u_i$, including the coordinate functions. Their product remains nonzero in the coordinates of an $L$-basis of $K_L$, and it cannot vanish on all of $L^6$, since $L$ is infinite. Let $t=u^6/N_{K_L/L}(u)$. Since $[K_L:L]=6$, we have $N_{K_L/L}(t)=1$ and $t_i/t_j=(u_i/u_j)^6$.

We claim that $W$ contains an element with a nonzero off-diagonal component. By Proposition~\ref{prop:base-change}(ii), we have $\Aut(K_L/L)=1$, hence the normalizer of $K_L$ in $D_L^\times$ is $C_{D_L}(K_L)^\times=K_L^\times$. Note that $t\notin L$ because the ratios $t_i/t_j$ with $i\ne j$ differ from $1$; hence, $L(t)=K_L$ since there are no proper intermediate fields. Let $h=xtx^{-1}\in H$. If $h\in K_L^\times$, then $xK_Lx^{-1}=L(h)=K_L$ by degree, contradicting the choice of $x\notin K_L^\times$. So, $h\notin K_L^\times$. Now, $K_L=L\cdot1\oplus K_L^0$ since $\tr_{K_L/L}(1)=6$. Hence, the inclusion $hK_L^0h^{-1}\subseteq K_L$ would imply $hK_Lh^{-1}=K_L$, which is impossible since $h\notin K_L^\times$. But $D_L\cap(K_L\otimes_L k)=K_L$ inside $A$. Therefore, some element of $hK_L^0h^{-1}\subseteq W$ must have a nonzero off-diagonal component, proving the claim.

Let $e_i$ be the $i$-th diagonal matrix unit, and let $c_t\colon A\to A$ be the map defined by $a\mapsto tat^{-1}$. Then, $c_t(W_k)=W_k$ since $t\in H$. Also, $c_t$ acts by $1$ on the diagonal matrices and by $t_i/t_j$ on $e_iAe_j$ for $i\ne j$. Interpolation at the distinct eigenvalues shows that each projection onto $e_iAe_j$ with $i\ne j$ belongs to $k[c_t]$, and therefore preserves $W_k$. By the claim, one such projection is nonzero on $hK_L^0h^{-1}$, hence $W_k$ contains its image line.

Consider now the semilinear action of $\Aut(k/L)$ on $A$. Every $\sigma\in\Aut(k/L)$ satisfies $\sigma(W_k)=W_k$ and $\sigma(e_iAe_j)=e_{\sigma(i)}Ae_{\sigma(j)}$, where the indices are permuted as the embeddings of $K_L$. By Proposition~\ref{prop:base-change}(ii), this is the $S_6$-action, which is transitive on ordered pairs of distinct indices. Therefore, $W_k$ contains every $e_iAe_j$ with $i\ne j$. It also contains the trace-zero diagonal matrices, which form $K_L^0\otimes_L k$, and it is contained in the kernel of the matrix trace. Hence,
\[
W_k=(K_L^0\otimes_L k)\oplus{\textstyle\bigoplus_{i\ne j}}\,e_iAe_j=\ker(\tr\colon A\to k),
\]
so $W=\ker\Trd_{D_L}$.
\end{proof}

In view of Lemma~\ref{lem:W}, finitely many conjugates of $T$ suffice to make the differential of the multiplication map surjective:

\begin{lemma}\label{lem:dense}
There are $h_1,\ldots,h_r\in H$, with $h_1=1$, for which the morphism
\begin{equation}\label{eq:phi}
\phi\colon T^r\longrightarrow V\quad\text{defined by}\quad
(t_1,\ldots,t_r)\longmapsto{\textstyle\prod_{j=1}^r}\,h_jt_jh_j^{-1}
\end{equation}
has surjective differential at $(1,\ldots,1)$. In particular, the morphism $\phi$ is dominant.
\end{lemma}

\begin{proof}
Since $W$ is finite dimensional over $L$, there are $h_1=1,h_2,\ldots,h_r\in H$ with $W=\sum_{j=1}^r h_jK_L^0h_j^{-1}$. Let $\phi$ be the corresponding morphism \eqref{eq:phi}.

Over $L[\varepsilon]/(\varepsilon^2)$, the norm polynomials satisfy
\[
N_{K_L/L}(1+\varepsilon a)=1+\varepsilon\tr_{K_L/L}(a)\quad\text{and}\quad
\Nrd_{D_L}(1+\varepsilon b)=1+\varepsilon\Trd_{D_L}(b).
\]
Hence, the tangent spaces at $1$ are $T_1T=K_L^0$ and $T_1V=\ker\Trd_{D_L}$. For $a_1,\ldots,a_r\in K_L^0$, multiplication in the definition of $\phi$ yields
\[
\phi(1+\varepsilon a_1,\ldots,1+\varepsilon a_r)
=1+\varepsilon{\textstyle\sum_{j=1}^r}\,h_j a_j h_j^{-1}.
\]
Therefore, the image of the differential $d\phi_{(1,\ldots,1)}$ is $W$, which equals $T_1V$ by Lemma~\ref{lem:W}.

It remains to show that $\phi$ is dominant. Let $k=L_{\alg}$. Over $k$, the defining equations of $T$ and $V$ are $t_1\cdots t_6=1$ and $\det(a)=1$. Both varieties are smooth, since the derivatives of the two polynomials in the directions $t$ and $a$ at their points are $6$. Also, $V$ is geometrically irreducible, because the morphism from the open set of invertible matrices to $V_k$ defined by $a\mapsto\operatorname{diag}(\det(a)^{-1},1,\ldots,1)a$ is surjective. By the differential criterion for smoothness, $\phi$ is smooth at $(1,\ldots,1)$. Its image therefore contains a nonempty open subset of $V$, hence $\phi$ is dominant.
\end{proof}

Henceforth, we fix elements $h_1=1,h_2,\ldots,h_r\in H$ and the morphism $\phi$ as in Lemma~\ref{lem:dense}.

\subsection{Local generation}\label{ss:local}

In this subsection, we show that the conjugates $h_jT(L_v)h_j^{-1}$ generate the full group of elements of reduced norm $1$ at each place $v$ of $L$ above $\Sigma$. The argument uses elementary matrices, for which we fix the following notation. Let $\Delta$ be a centrally finite division ring, and let $m\ge2$. For $i\ne j$ and $b\in\Delta$, let $u_{ij}(b)=I+bE_{ij}\in M_m(\Delta)$, and let
\begin{equation}\label{eq:Em}
E_m(\Delta)=\langle u_{ij}(b)\mid i\ne j,\ b\in\Delta\rangle.
\end{equation}
Over a splitting field of $\Delta$, the matrix $u_{ij}(b)$ is block triangular with identity diagonal blocks, hence $\Nrd_{M_m(\Delta)}(u_{ij}(b))=1$. This shows that $E_m(\Delta)\subseteq M_m(\Delta)^{(1)}$. The next lemma compares these two groups by means of the Dieudonn\'e determinant $\det$ (see Tignol--Wadsworth \cite[Lemma~11.22, p.~549]{TignolWadsworth2015}).

\begin{lemma}\label{lem:elementary}
Let $\Delta$ be a centrally finite division ring, and let $m\ge2$. Then, the Dieudonn\'e determinant induces an isomorphism
\begin{equation}\label{eq:SK1}
M_m(\Delta)^{(1)}/E_m(\Delta)\xrightarrow[\sim]{\ \det\ }\SK_1(\Delta).
\end{equation}
In particular, if $\SK_1(\Delta)=1$, then $M_m(\Delta)^{(1)}=E_m(\Delta)$.
\end{lemma}

\begin{proof}
Since $m\ge2$, the Dieudonn\'e determinant fits into the exact sequence of Dieudonn\'e \cite[Theorem~1, pp.~32--33]{Dieudonne}
\begin{equation}\label{eq:dieudonne}
1\longrightarrow E_m(\Delta)\longrightarrow M_m(\Delta)^\times\xrightarrow{\ \det\ }\Delta^\times/[\Delta^\times,\Delta^\times]\longrightarrow 1.
\end{equation}
Let $\overline{\Nrd}_\Delta$ denote the map induced by $\Nrd_\Delta$ on the abelianization of $\Delta^\times$. Both $\Nrd_{M_m(\Delta)}$ and $\overline{\Nrd}_\Delta\circ\det$ are trivial on $E_m(\Delta)$. For $d(c)=\operatorname{diag}(c,1,\ldots,1)$ with $c\in\Delta^\times$, we have
\[
\Nrd_{M_m(\Delta)}(d(c))=\Nrd_\Delta(c)=\overline{\Nrd}_\Delta(\det d(c))
\]
by Draxl \cite[\S22, Theorem~1 and Corollary~1, pp.~146--148]{Draxl1983}. But the matrices $d(c)$ represent all cosets modulo $E_m(\Delta)$ by \eqref{eq:dieudonne}, hence $\Nrd_{M_m(\Delta)}=\overline{\Nrd}_\Delta\circ\det$. Therefore, $\det$ maps $M_m(\Delta)^{(1)}$ onto $\ker\overline{\Nrd}_\Delta=\SK_1(\Delta)$, with kernel $E_m(\Delta)$. This proves \eqref{eq:SK1}.
\end{proof}

We now prove the local generation statement.

\begin{lemma}\label{lem:local}
For each place $v$ of $L$ above $\Sigma$, we have
\begin{equation}\label{eq:local}
(D_L\otimes_L L_v)^{(1)}=\langle h_jT(L_v)h_j^{-1}\mid 1\le j\le r\rangle.
\end{equation}
\end{lemma}

\begin{proof}
Let $k=L_v$, and let $J$ be the group on the right side of \eqref{eq:local}. We proceed in three steps.

\emph{Step 1: $J$ contains an open neighborhood of $1$ in $(D_L\otimes_L k)^{(1)}$.} We have $\phi(T(k)^r)\subseteq J$ by the definition \eqref{eq:phi} of $\phi$. Since the differential of $\phi$ at $(1,\ldots,1)$ is surjective over $k$ by Lemma~\ref{lem:dense}, the implicit function theorem over $k$ shows that $\phi(T(k)^r)$ contains an open neighborhood of $1$ in $(D_L\otimes_L k)^{(1)}$.

\emph{Step 2: $J$ contains the elementary matrices, for a suitable identification of $D_L\otimes_L k$ with a matrix algebra.} Let $p$ be the prime below $v$. By Proposition~\ref{prop:base-change}(iii), we may identify $k$ with $\Q_p$, and we have
\[
D_L\otimes_L k\simeq M_m(\Delta)\quad\text{and}\quad K_L\otimes_L k\simeq{\textstyle\prod_{\ell=1}^m}\,P_{p,\ell},
\]
where $m=m_p\ge2$, $\Delta=\Delta_p$, and $[P_{p,\ell}:k]=d_p$. The idempotents $e_\ell$ of the field factors yield a decomposition $\Delta^m=\bigoplus_{\ell=1}^m e_\ell\Delta^m$. Since all $m$ summands are nonzero, each has dimension one over $\Delta$. We fix a basis of $\Delta^m$ adapted to these summands. Take any $i\ne j$, let $\pi$ be a uniformizer of $k$, and let
\[
a=\pi e_i+\pi^{-1}e_j+\sum_{\ell\ne i,j}e_\ell.
\]
This matrix is diagonal, and its norm as an element of $K_L\otimes_L k$ is $\pi^{d_p}\pi^{-d_p}=1$. Hence, $a\in T(k)\subseteq J$, since $h_1=1$. Let $b\in\Delta$. Since $\pi$ is central, we have
\[
a^Nu_{ij}(b)a^{-N}=I+(\pi^Nb\pi^N)E_{ij}
=u_{ij}(\pi^{2N}b)\longrightarrow I\quad\text{as }N\to\infty.
\]
Since $J$ contains a neighborhood of $I$ by Step~1, the conjugate $a^Nu_{ij}(b)a^{-N}$ lies in $J$ for large $N$. As $a\in J$, it follows that $u_{ij}(b)=a^{-N}(a^Nu_{ij}(b)a^{-N})a^N\in J$. This shows that $E_m(\Delta)\subseteq J$.

\emph{Step 3: $J=(D_L\otimes_L k)^{(1)}$.} Since $\ind(\Delta)=d_p\in\{2,3\}$ is squarefree, we have $\SK_1(\Delta)=1$ by Wang's Theorem (see Tignol--Wadsworth \cite[Proposition~11.2, p.~536]{TignolWadsworth2015}). Therefore, Lemma~\ref{lem:elementary} and Step~2 yield
\[
E_m(\Delta)\subseteq J\subseteq M_m(\Delta)^{(1)}=E_m(\Delta),
\]
which completes the proof.
\end{proof}

\subsection{Fibers of convolution powers}\label{ss:fibers}

Let $Y=T^r$, and let $e=(1,\ldots,1)\in Y(L)$. For $q\ge1$, let $\phi^{*q}\colon Y^q\to V$ be the morphism defined by
\begin{equation}\label{eq:conv}
\phi^{*q}(y_1,\ldots,y_q)=\phi(y_1)\cdots\phi(y_q).
\end{equation}
Factoring an element $s\in D_L^{(1)}$ as in Proposition~\ref{prop:factor} amounts to finding a point on the fiber
\begin{equation}\label{eq:Xs}
X_s=(\phi^{*q})^{-1}(s)
\end{equation}
of such a map over $s$, with coordinates in a suitable extension of $L$. The goal of this subsection is to prove that these fibers remain irreducible over $L_{\alg}=\Q_{\alg}$ when $q$ is large. The proof relies on results of Breuillard--Glazer \cite{BG} on convolutions of morphisms to algebraic groups over $\Q$, which we apply after restricting scalars from $L$ to $\Q$: in a $\Q$-basis of $L$, each equation over $L$ becomes $[L:\Q]$ equations over $\Q$ by equating coefficients. This is the Weil restriction $R_{L/\Q}$ of Tignol--Wadsworth \cite[\S12.3, p.~577]{TignolWadsworth2015}.

\begin{lemma}\label{lem:fibers}
Let $s\in D_L^{(1)}$, and let $q\ge35[L:\Q]+1$. If the fiber $X_s$ of \eqref{eq:Xs} is nonempty, then it remains irreducible over $\Q_{\alg}$.
\end{lemma}

\begin{proof}
Let
\[
G=R_{L/\Q}(V)\quad\text{and}\quad
\psi=R_{L/\Q}(\phi)\colon R_{L/\Q}(Y)\longrightarrow G.
\]
Let $k=\Q_{\alg}$, and let $\sigma$ range over the embeddings $L\hookrightarrow k$. The decomposition $L\otimes_\Q k\simeq{\textstyle\prod_\sigma}\,k$ identifies $\psi_k$ with the product of the scalar extensions $\phi_\sigma=\phi\times_{L,\sigma}k$.

Each source $Y\times_{L,\sigma}k$ is isomorphic to $(\mathbb A^1_k\setminus\{0\})^{5r}$. Hence, $R_{L/\Q}(Y)$ is smooth and geometrically irreducible, and $\psi$ is dominant by Lemma~\ref{lem:dense} applied to each factor. Moreover, $G$ is an algebraic group: multiplication on $G$ is induced by that of $D_L$, and inversion on $V$ is regular by the reduced characteristic polynomial and the equation $\Nrd_{D_L}=1$. We have
\[
G_k\simeq{\textstyle\prod_\sigma}\,\SL_6\quad\text{and}\quad\dim G=35[L:\Q].
\]
In particular, $G$ is connected and absolutely simply connected. Let $e_0\in R_{L/\Q}(Y)(\Q)$ correspond to $e$, so $\psi(e_0)=1$. Then, $\psi^{*t}(y,e_0,\ldots,e_0)=\psi(y)$ for $t\ge1$, which shows that every $\psi^{*t}$ is dominant.

Since $q-1\ge\dim G\ge2$, the morphism $\psi^{*(q-1)}$ is flat by Breuillard--Glazer \cite[Theorem~4.16, p.~101]{BG}. Since $q-2\ge1$, its generic fiber is geometrically irreducible by \cite[Theorem~8.2, p.~123]{BG}, applied to $\psi^{*(q-1)}=\psi*\psi^{*(q-2)}$. It follows by \cite[Proposition~4.21, p.~102]{BG} that every fiber of $\psi^{*q}=\psi^{*(q-1)}*\psi$ is geometrically irreducible.

Finally, let $s_0\in G(\Q)$ correspond to $s$. The decomposition of $L\otimes_\Q k$ above identifies the geometric fiber of $\psi^{*q}$ over $s_0$ with
\[
\bigl((\psi^{*q})^{-1}(s_0)\bigr)_k
\simeq{\textstyle\prod_\sigma}\,(X_s\times_{L,\sigma}k).
\]
Since $X_s$ is nonempty, each factor has a $k$-point. Hence, the projection onto the factor for the inclusion $L\subseteq k$ is surjective. But the product is irreducible, so $(X_s)_k$ must be irreducible.
\end{proof}

\subsection{Proof of \texorpdfstring{Proposition~\ref{prop:factor}}{Proposition 4.1}}\label{ss:proof-factor}

Let $s\in D_L^{(1)}$. We first factor $s$ locally. For every place $v$ of $L$ above $\Sigma$, the image $\phi(Y(L_v))$ contains each $h_jT(L_v)h_j^{-1}$, as we see by taking all coordinates except the $j$-th equal to $1$. By Lemma~\ref{lem:local}, we may therefore write
\begin{equation}\label{eq:local-fact}
s=\phi(y_{v,1})\cdots\phi(y_{v,b_v}),\quad\text{where }y_{v,i}\in Y(L_v).
\end{equation}
Choose an integer
\begin{equation}\label{eq:q}
q\ge\max\{35[L:\Q]+1,\ 1+\max_{v\mid\Sigma}b_v\},
\end{equation}
and let $X_s=(\phi^{*q})^{-1}(s)$, as in \eqref{eq:Xs}. Then, \eqref{eq:local-fact} shows that, for each place $v$ of $L$ above $\Sigma$,
\begin{equation}\label{eq:zv}
z_v=(y_{v,1},\ldots,y_{v,b_v},e,\ldots,e)\in X_s(L_v).
\end{equation}

The point $z_v$ is smooth on $(X_s)_{L_v}$. To see this, let $\lambda_s\colon V_{L_v}\to V_{L_v}$ be the automorphism defined by $a\mapsto sa$. By \eqref{eq:q}, the last coordinate of $z_v$ is $e$. Keeping the earlier factors at their values in $z_v$ restricts $\phi^{*q}$ over $L_v$ to $\lambda_s\circ\phi$ on the last factor. Hence,
\[
\left.d(\phi^{*q})_{z_v}\right|_{T_e(Y_{L_v})}=d(\lambda_s)_1\circ d\phi_e.
\]
The right side is surjective, since $d(\lambda_s)_1$ is an isomorphism and $d\phi_e$ is surjective by Lemma~\ref{lem:dense}. Since $Y^q$ and $V$ are smooth, the differential criterion yields the smoothness of $(X_s)_{L_v}$ at $z_v$.

By \eqref{eq:zv}, the fiber $X_s$ is nonempty, hence it remains irreducible over $\Q_{\alg}$ by Lemma~\ref{lem:fibers}, in view of \eqref{eq:q}. So, the reduced fiber $X=(X_s)_{\mathrm{red}}$ is an affine variety over $L$ that remains irreducible over $L_{\alg}$. Moreover, a smooth neighborhood of $z_v$ in $(X_s)_{L_v}$ is also open in $X_{L_v}$, so $z_v\in X_{\mathrm{sm}}(L_v)$ for each place $v$ of $L$ above $\Sigma$. By Lemma~\ref{lem:MB}, there is therefore a finite extension $L'/L$ inside $F$ with $X(L')\ne\varnothing$. Since $T(L)\subseteq T(L')$, the subgroup $\langle T(L'),xT(L')x^{-1}\rangle$ contains $H$, hence all the $h_j$. Thus, for $(y_1,\ldots,y_q)\in X(L')$, we have
\[
s=\phi(y_1)\cdots\phi(y_q)\in\langle T(L'),xT(L')x^{-1}\rangle,
\]
which completes the proof of Proposition~\ref{prop:factor}. \qed

\section{Proof of \texorpdfstring{Theorem~\ref{thm:main}}{Theorem 1.1} and complements}\label{sec:proof}

In this section, we derive Theorem~\ref{thm:main} from Proposition~\ref{prop:norm} and Proposition~\ref{prop:factor}. We then determine in \S\ref{ss:structure} the normal subgroups of $D^{(1)}$, and we observe in \S\ref{ss:valuations} that $D$ carries no nontrivial valuation. Throughout, $D=D_F$ and $K=K_F$ are as in \S\ref{ss:scalar}.

\subsection{Proof of \texorpdfstring{Theorem~\ref{thm:main}}{Theorem 1.1}}\label{ss:proof-main}

(i) By Proposition~\ref{prop:base-change} and Proposition~\ref{prop:cyclic}, $D$ is a cyclic division algebra of degree $6$ with center $F$, and $K$ is a maximal subfield whose normal closure over $F$ has Galois group $S_6$. Moreover, $N_{K/F}$ is surjective by Proposition~\ref{prop:norm}.

(ii) Let $x\in D^\times\setminus K^\times$. Since $\ker N_{K/F}=K^\times\cap D^{(1)}$ by \eqref{eq:kerN} and $D^{(1)}\trianglelefteq D^\times$, the subgroup on the left in (ii) is contained in $D^{(1)}$. For the reverse inclusion, take any $s\in D^{(1)}$, and let $L\subseteq F$ be a number field with $s,x\in D_L$. Then, $s\in D_L^{(1)}$ by compatibility of the reduced norm with scalar extension, and $x\in D_L^\times\setminus K_L^\times$ since $K_L\subseteq K$. Hence, by Proposition~\ref{prop:factor} there is a finite extension $L'/L$ inside $F$ such that
\[
s\in\langle T(L'),xT(L')x^{-1}\rangle
\subseteq\langle\ker N_{K/F},x(\ker N_{K/F})x^{-1}\rangle,
\]
where the inclusion follows from $T(L')=\ker N_{K_{L'}/L'}\subseteq\ker N_{K/F}$, by compatibility of field norms with scalar extension. This proves (ii).

(iii) Let $y\in D^\times$. By (i), there is some $a\in K^\times$ with $N_{K/F}(a)=\Nrd_D(y)$. Then, $a^{-1}y\in D^{(1)}$. Hence,
\begin{equation}\label{eq:product}
D^\times=K^\times D^{(1)}.
\end{equation}
Note that $K^\times\subsetneq D^\times$ since $[K:F]=6<36=[D:F]$. In view of \eqref{eq:product}, we therefore have $\ker N_{K/F}=K^\times\cap D^{(1)}\subsetneq D^{(1)}$.

If $x\in D^\times\setminus K^\times$, then $D^{(1)}\subseteq\langle K^\times,x\rangle$ by (ii); hence, $\langle K^\times,x\rangle=D^\times$ by \eqref{eq:product}. Likewise, if $x\in D^{(1)}\setminus\ker N_{K/F}$, then $x\notin K^\times$, so $\langle\ker N_{K/F},x\rangle=D^{(1)}$ by (ii). This proves both maximality assertions. \qed

\subsection{The group \texorpdfstring{$D^{(1)}/Z(D^{(1)})$}{D(1)/Z(D(1))}}\label{ss:structure}

The generation formula in Theorem~\ref{thm:main}(ii) relates the abelianization of $D^{(1)}$ to $\ker N_{K/F}$. Together with the malnormality of $K^\times/F^\times$, it yields the following information on the structure of $D^{(1)}$:

\begin{corollary}\label{cor:structure}
Let $D$ and $K$ be as in \S\ref{ss:scalar}.
\begin{enumerate}[label=(\roman*)]
\item The subgroup $K^\times/F^\times$ is malnormal in $D^\times/F^\times$.
\item We have $\SK_1(D)=1$, and the group $D^{(1)}$ is perfect.
\item The center of $D^{(1)}$ is $Z(D^{(1)})=\{\pm1\}$.
\item The quotient $D^{(1)}/Z(D^{(1)})$ is nonabelian simple, and the image of $\ker N_{K/F}$ is an abelian maximal subgroup.
\end{enumerate}
\end{corollary}

\begin{proof}
(i) Suppose $x\in D^\times\setminus K^\times$ and $a\in(K^\times\cap xK^\times x^{-1})\setminus F^\times$. Then, $F(a)=K=xKx^{-1}$ by Proposition~\ref{prop:base-change}(ii). But $\Aut(K/F)=1$ by the same proposition, so $x\in C_D(K)^\times=K^\times$, a contradiction. Hence,
\begin{equation}\label{eq:malnormal}
K^\times\cap xK^\times x^{-1}=F^\times\quad\text{for every }x\in D^\times\setminus K^\times,
\end{equation}
proving (i).

(ii) Since $\ind D=6$ is squarefree, the vanishing theorem of Wang \cite{Wang} applies (see also Tignol--Wadsworth \cite[Proposition~11.2, p.~536]{TignolWadsworth2015}). Thus, $\SK_1(D)=1$ and $D^{(1)}=[D^\times,D^\times]$. Let $C=[D^{(1)},D^{(1)}]$, which is a characteristic subgroup of $D^{(1)}$, hence a normal subgroup of $D^\times$. By Theorem~\ref{thm:main}(iii) and \eqref{eq:kerN}, there is some $x\in D^{(1)}\setminus K^\times$. Conjugation by $x$ acts trivially on $D^{(1)}/C$, hence $D^{(1)}=(\ker N_{K/F})C$ by Theorem~\ref{thm:main}(ii). By \eqref{eq:product}, we then have $D^\times=K^\times D^{(1)}=K^\times C$. Therefore, $D^\times/C$ is abelian and $D^{(1)}=[D^\times,D^\times]\subseteq C$, i.e., $D^{(1)}$ is perfect.

(iii) Let $N$ be a proper normal subgroup of $D^{(1)}$. Suppose $N\not\subseteq\ker N_{K/F}$. Then, $N(\ker N_{K/F})=D^{(1)}$ by the maximality in Theorem~\ref{thm:main}(iii), hence $D^{(1)}/N$ is nontrivial and abelian. This contradicts (ii). Therefore, $N\subseteq\ker N_{K/F}$. With $x\in D^{(1)}\setminus K^\times$ as above, normality and \eqref{eq:malnormal} imply
\begin{equation}\label{eq:normal}
N=xNx^{-1}\subseteq K^\times\cap xK^\times x^{-1}=F^\times.
\end{equation}
Since $D^{(1)}$ is nontrivial and perfect, its center is a proper subgroup. Applying \eqref{eq:normal} to $Z(D^{(1)})$, we obtain
\[
Z(D^{(1)})=F^\times\cap D^{(1)}=\{c\in F^\times\mid c^6=1\}\subseteq\ker N_{K/F},
\]
because $\Nrd_D(c)=c^6$ for $c\in F^\times$. Let $c\in F^\times$ with $c^6=1$. By Lemma~\ref{lem:center}(i), we may embed $\Q(c)$ in $\Q_2$; let $z$ be the image of $c^2$. Since $z^3=1$, we have $z\in\mathbb Z_2^\times$, hence $z\equiv1\pmod 2$. Thus,
\[
0=z^3-1=(z-1)(z^2+z+1),\quad\text{with}\quad z^2+z+1\equiv1\pmod 2.
\]
It follows that $z=1$, so $c^2=1$. Since $(\pm1)^6=1$, we have $Z(D^{(1)})=\{\pm1\}$, proving (iii).

(iv) By \eqref{eq:normal}, every proper normal subgroup of $D^{(1)}$ lies in $F^\times\cap D^{(1)}=Z(D^{(1)})$. Therefore, $D^{(1)}/Z(D^{(1)})$ is simple; it is nonabelian because it is nontrivial and perfect. Furthermore, the image of the abelian subgroup $\ker N_{K/F}$ is maximal in this quotient by Theorem~\ref{thm:main}(iii) and the subgroup correspondence.
\end{proof}

\subsection{Valuations on \texorpdfstring{$D$}{D}}\label{ss:valuations}

We conclude by considering valuations on $D$, in the sense of \S\ref{ss:K0}. By Proposition~\ref{prop:base-division}, the algebra $D_0\otimes_\Q\Q_p$ is not a division ring for any prime $p$. This prevents $D$ from carrying a nontrivial valuation:

\begin{proposition}\label{prop:valuations}
Every valuation on $D$ is trivial.
\end{proposition}

\begin{proof}
Let $w$ be a valuation on $D$, and let $a\in D^\times$. Then, $a$ is algebraic over $\Q$, since it lies in some $D_L$ with $[L:\Q]<\infty$. In a relation $\sum_i c_i a^i=0$ over $\Q$, the least value of a nonzero term must occur at least twice. Thus,
\begin{equation}\label{eq:val}
(j-i)w(a)=w(c_i)-w(c_j)\quad\text{for some }i<j\text{ with }c_i c_j\ne0.
\end{equation}
In particular, if $w$ is trivial on $\Q$, then it is trivial on $D$.

Suppose now that $w$ is nontrivial. Then, $w$ is nontrivial on $\Q$ and nonnegative on nonzero integers; so, prime factorization and B\'ezout's identity show that there is a unique prime $p$ with $w(p)>0$. Since $w(q)=v_p(q)w(p)$ for $q\in\Q^\times$, it follows from \eqref{eq:val} that $(j-i)w(a)\in\mathbb Z w(p)$. We may thus identify the value group of $w$ with a subgroup of $\Q$ via the embedding that sends $w(p)$ to $1$.

Let $\widehat{D_0}$ be the completion of $D_0$ for $w|_{D_0}$. We claim that $\widehat{D_0}$ is a division ring. Let $(a_i)$ be a sequence representing a nonzero element of $\widehat{D_0}$. Then, $w(a_i)=\gamma$ eventually, and
\[
w(a_i^{-1}-a_j^{-1})=w\bigl(a_i^{-1}(a_j-a_i)a_j^{-1}\bigr)
=w(a_j-a_i)-2\gamma\longrightarrow\infty.
\]
Hence, the inverses form a Cauchy sequence, proving the claim. Now, the closure of $\Q$ in $\widehat{D_0}$ is a central copy of $\Q_p$. Since $D_0\otimes_\Q\Q_p$ is simple, multiplication induces a unital embedding
\[
D_0\otimes_\Q\Q_p\hookrightarrow\widehat{D_0}\quad\text{defined by}\quad d\otimes q\longmapsto dq.
\]
But $D_0\otimes_\Q\Q_p$ has matrix size $m_p>1$ for $p\in\Sigma$ and $6$ otherwise, by Proposition~\ref{prop:base-division}. So, a nontrivial idempotent would embed in the division ring $\widehat{D_0}$, a contradiction. Therefore, $w$ must be trivial.
\end{proof}

\end{document}